\documentclass[11pt,reqno]{amsart}
\usepackage{amsmath}
\usepackage{amssymb, latexsym,
amsthm,
amsfonts, amsbsy,amsthm, amsmath, mathrsfs}
\usepackage{bm}
\usepackage{graphicx}

\newtheorem{thm}{Theorem}

\newtheorem{lem}[thm]{Lemma}
\newtheorem{cor}[thm]{Corollary}
\newtheorem{remk}[thm]{Remark}

\newtheorem{defn}[thm]{Definition}

\newtheorem{exmp}[thm]{Example}

\newcommand{\de}{\delta}

\newcommand{\ep}{\varepsilon}
\newcommand{\rig}{\rightarrow}
\newcommand{\mrig}{\mathrel{-\!\!\!\!\!\rightarrow}}

\newcommand{\sbA}{{\boldsymbol{A}}}
\newcommand{\sbB}{{\boldsymbol{B}}}

\newcommand{\sbT}{{\boldsymbol{T}}}

\newcommand{\Eq}{{\boldsymbol{\mathit{Eq}}}}
\newcommand{\Con}{{\boldsymbol{\mathit{Con}}}}
\newcommand{\Rel}{{\boldsymbol{\mathit{Rel}}}}

\newcommand{\ov}{\overline}

\let\class=\mathsf                              
\bmdefine{\A}{A}                                

\newcommand{\HHH}{\mathbb{H}}

\newcommand{\PPP}{\mathbb{P}}

\newcommand{\PPU}{\mathbb{P}_{\!\textsc{u}}^{}}

\newcommand{\SSS}{\mathbb{S}}
\newcommand{\III}{\mathbb{I}}

\bmdefine{\btau}{\tau}                                  
\bmdefine{\brho}{\rho}                                  

\newcommand{\semeq}{\mathrel{=\joinrel\mathrel\vert\mkern0mu\mathrel\vert\joinrel=}}

\newcommand{\seteq}{\mathrel{\mbox{\,\textup{:}\!}=\nolinebreak }\,}
\newcommand{\coloneqq}{\mathrel{\mbox{\,\textup{:}\!}=\nolinebreak }\,}

\theoremstyle{theorem}
\newtheorem{Theorem}{Theorem}[section]
\newtheorem{Lemma}[Theorem]{Lemma}

\theoremstyle{definition}

\newtheorem{exa}[Theorem]{Example}

\theoremstyle{remark}

\begin{document}
\title[Prevarieties of Logic]{On Prevarieties of Logic} 

\author{Tommaso Moraschini}
\address{Institute of Computer Science, Academy of Sciences of the Czech Republic, Pod Vod\'{a}renskou v\v{e}\v{z}\'{i} 2, 182 07 Prague 8, Czech Republic.}
\email{moraschini@cs.cas.cz}
\author{James G.\ Raftery}
\address{Department of Mathematics and Applied Mathematics,
 University of Pretoria,
 Private Bag X20, Hatfield,
 Pretoria 0028, South Africa}
\email{james.raftery@up.ac.za}


\thanks{This work received funding from the European Union's Horizon 2020 research and innovation programme under the Marie Sklodowska-Curie grant
agreement No~689176 (project ``Syntax Meets Semantics: Methods, Interactions, and Connections in Substructural logics").
The first author was supported by project CZ.02.2.69/0.0/0.0/17\_050/0008361, OPVVV M\v{S}MT, MSCA-IF Lidsk\'{e} zdroje v teoretick\'{e} informatice.
The second author was supported in part by the National Research Foundation (NRF) of South Africa (UID 85407).
Both authors thank the University of Pretoria and
the DST-NRF Centre of Excellence in Mathematical and Statistical Sciences (CoE-MaSS)
for
partially funding the first author's travel to Pretoria in 2017 and 2018}


\begin{abstract}
It is proved that every prevariety of algebras is categorically equivalent to a `prevariety of logic', i.e., to the
equivalent algebraic semantics of some sentential deductive system.
This allows us to show that no nontrivial
equation in the language $\land,\lor,\circ$ holds in the congruence lattices of all members of every
variety of logic, and that being a (pre)variety of logic is not a categorical property.
\end{abstract}

$\mbox{}$\\{\vspace{-8mm}}

\maketitle

\makeatletter
\renewcommand{\labelenumi}{\text{(\theenumi)}}
\renewcommand{\theenumi}{\roman{enumi}}
\renewcommand{\theenumii}{\roman{enumii}}
\renewcommand{\labelenumii}{\text{(\theenumii)}}
\renewcommand{\p@enumii}{\theenumi(\theenumii)}
\makeatother

\allowdisplaybreaks

\section{Prevarieties of Logic}

Recall that the class operator symbols $\III$, $\HHH$, $\SSS$, $\PPP$ and $\PPU$ stand for the formation of isomorphic and homomorphic images,
subalgebras, direct products and ultraproducts, respectively.  A class of similar algebras is called a \emph{pre\-variety}, a
\emph{quasivariety} or a \emph{variety} if it is closed, respectively, under $\III$, $\SSS$ and $\PPP$, under
$\III$, $\SSS$, $\PPP$ and $\PPU$, or under $\HHH$, $\SSS$ and $\PPP$.

The informal notion of a `variety of logic' has acquired a precise meaning in abstract algebraic logic (see \cite{BP89,Cze01,Fon16}),
where it extends naturally
to prevarieties.  In the standard terminology, a
prevariety of logic
is the equivalent algebraic semantics
of an algebraizable (sentential) logic, but the following purely algebraic characterization can serve here as a definition.
\begin{defn}\label{def:prevariety of logic}
\textup{A prevariety $\class{K}$ is called a \emph{prevariety of logic} if some fixed formula of infinitary logic, having the form
\begin{equation}\label{eq:prevariety of logic}
\big(\text{\large \&}_{i\in I,\,j\in J}\,\,\de_i(\rho_j(x,y))\approx\ep_i(\rho_j(x,y))\big)\,\Longleftrightarrow \,x\approx y,
\end{equation}
is valid in (every member of) $\class{K}$.  It is understood here that $I$ and $J$ are sets, and that
$\btau=\{\langle\delta_i,\ep_i\rangle:i\in I\}$ is a family of pairs of unary terms
and $\brho=\{\rho_j:j\in J\}$ a family of binary terms in the signature of $\class{K}$.
In this context, $\btau$ and $\brho$ are called \emph{transformers}.  If, moreover, $\class{K}$ is a [quasi]variety, then we
refer to it as a [\emph{quasi}]\emph{variety of logic}.}\qed
\end{defn}
\begin{exmp}\label{boolean etc}
\textup{In
the variety of Boolean [resp.\ Heyting] algebras, which algebraizes classical [resp.\ intuitionistic] propositional logic,
(\ref{eq:prevariety of logic})
takes the form
\[
\left(x\rig y\approx 1 \;\;\&\;\; y\rig x\approx 1\right)\Longleftrightarrow \,x\approx y.
\]
In the variety of commutative residuated lattices \cite{GJKO07}, which algebraizes a rich fragment of linear logic, (1) is most
naturally instantiated as}

\smallskip

\quad\quad\quad\quad $\left((x\rig y)\wedge 1\approx 1 \;\;\&\;\; (y\rig x)\wedge 1\approx 1\right)\Longleftrightarrow \,x\approx y$. \qed
\end{exmp}
\indent
Given $\class{K},\btau$ and $\brho$ as in Definition~\ref{def:prevariety of logic}, we can construct a logic $\,\vdash_{\class{K},\btau}$
for which $\class{K}$ is the equivalent algebraic semantics, as follows.
It is convenient here to fix
a proper class $\mathit{Var}$ of
variables
for the entire discussion.

For each set $X\subseteq\mathit{Var}$,
a term $\varphi$ over
$X$ in the signature of $\class{K}$
is declared a $\,\vdash^X_{\class{K},\btau}$-\emph{consequence} of a set $\Gamma$ of such terms
(written as $\Gamma\vdash^X_{\class{K},\btau}\varphi$)
provided that the following is true:
for any homomorphism $h$ from the absolutely free algebra $\sbT(X)$ over $X$ to any member of $\class{K}$, the kernel of $h$ contains
\[
\btau(\varphi)\seteq\{\langle\de_i(\varphi),\ep_i(\varphi)\rangle:i\in I\}
\]
whenever it contains $\btau[\Gamma]\seteq\bigcup_{\gamma\in\Gamma}\btau(\gamma)$.
(This criterion is abbreviated
as
\begin{equation}\label{tau interpretation}
\btau[\Gamma]\models_\class{K}\btau(\varphi).)
\end{equation}
Thus, $\,\vdash^X_{\class{K},\btau}$ is a binary relation from the power set of $T(X)$ to $T(X)$.

For any two sets $X,Y\subseteq\mathit{Var}$, with $\Gamma\cup\{\varphi\}\subseteq T(X)\cap T(Y)$,
it can be verified that $\Gamma\vdash^X_{\class{K},\btau}\varphi$ iff $\Gamma\vdash^Y_{\class{K},\btau}\varphi$.
It therefore makes sense to write
\[
\textup{$\Gamma\vdash_{\class{K},\btau}\varphi$ if there
exists a set $X\subseteq\mathit{Var}$ such that $\Gamma\vdash^X_{\class{K},\btau}\varphi$.}
\]
Technically, $\,\vdash_{\class{K},\btau}$ is the family of relations $\,\vdash^X_{\class{K},\btau}$
indexed by the subsets $X$ of $\mathit{Var}$.
It has the following properties
for any sets $X,Y\subseteq\mathit{Var}$, any
$\Gamma\cup\Psi\cup\{\varphi\}\subseteq T(X)$ and any homomorphism $h\colon\sbT(X)\mrig\sbT(Y)$,
where we abbreviate
$\,\vdash_{\class{K},\btau}$ as $\,\vdash$\,:
\begin{enumerate}
\item\label{reflexivity}
if $\varphi\in\Gamma$, then $\Gamma\vdash\varphi$;

\smallskip

\item\label{transitivity}
if $\Gamma\vdash\psi$ for all $\psi\in\Psi$, and $\Psi\vdash\varphi$, then $\Gamma\vdash\varphi$;

\smallskip

\item\label{substitution invariance2}
if $\Gamma
\vdash
\varphi$, then $h[\Gamma]
\vdash
h(\varphi)$.
\setcounter{newexmp}{\value{enumi}}
\end{enumerate}

For present purposes, (\ref{reflexivity})--(\ref{substitution invariance2}) are the defining properties of \emph{logics}
(over $\mathit{Var}$) in general.  Notice that $\,\vdash_{\class{K},\btau}$ is defined and is a logic for \emph{any}
class $\class{K}$ of similar algebras and \emph{any} set $\btau$ of pairs of unary terms in its signature (regardless of
$\brho$ and (\ref{eq:prevariety of logic})).

We say that a logic $\,\vdash$ is \emph{finitary} if it has the following additional property:
\begin{enumerate}
\setcounter{enumi}{\value{newexmp}}
\item\label{finitarity}
whenever $\Gamma\vdash\varphi$, then $\Gamma'\vdash\varphi$ for some finite $\Gamma'\subseteq\Gamma$.
\end{enumerate}
In this case, for any infinite set $X\subseteq\mathit{Var}$, the logic $\,\vdash$ is determined by its restriction
to terms over $X$,
and (\ref{substitution invariance2}) need only be
stipulated for the endomorphisms $h$ of $\sbT(X)$.  A finitary logic is usually (and can always be) specified by a
formal system $\mathbf{F}$ of axioms and finite inference rules, where the natural deducibility relation of $\mathbf{F}$
satisfies (\ref{reflexivity})--(\ref{finitarity}).
\begin{defn}
\textup{If a logic $\,\vdash$ has the form $\,\vdash_{\class{K},\btau}$ for some prevariety $\class{K}$ and
transformers $\btau$ and $\brho$, where $\class{K}$ satisfies (\ref{eq:prevariety of logic}), then $\class{K}$
is said to \emph{algebraize} $\,\vdash$, and $\,\vdash$ is said to be \emph{algebraizable}.\,\footnote{\,This
implies that $\,\vdash$ and the full equational consequence relation $\models_{\,\class{K}}$ of $\class{K}$ are essentially
interchangeable (i.e., the interpretation given by (\ref{tau interpretation}) is invertible), but we shall not
dwell further on the interpretations here.  See \cite{Fon16} for additional motivation.}}\qed
\end{defn}
\noindent
Under these conditions, $\class{K}$ is called the
\emph{equivalent algebraic semantics} of $\,\vdash$, because it is the \emph{only} prevariety that algebraizes $\,\vdash$;
the transformers are essentially unique as well,
cf.\ \cite[Thm.~2.15]{BP89}.  A quasivariety that algebraizes two
different logics (via different transformers $\btau$ and a common $\brho$) is exhibited in \cite[Sec.~5.2]{BP89},
along with several non-algebraizable finitary logics.

An algebraizable logic is said to be \emph{finitely algebraizable} if, in its algebraization, the transformer $\brho$
can be chosen finite,
i.e., the index set
$J$ in (\ref{eq:prevariety of logic}) can be kept
finite.  That happens, for instance, whenever the equivalent algebraic semantics $\class{K}$ is a quasivariety.  Dually, if an
algebraizable logic $\,\vdash$
is finitary,
then a finite choice of $\btau$ (i.e., of $I$) is possible in (\ref{eq:prevariety of logic}).
These facts and the next lemma are proved, for instance, in
\cite[Lem.~3.37]{Fon16}.
\begin{lem}
\label{Lem:algebraizability}
Let\/ $\class{K}$ be a quasivariety of logic, with transformers\/ $\btau$ and\/ $\brho$
as in \textup{(\ref{eq:prevariety of logic}).}
Then\/ $\btau$ can be chosen finite
iff\/ $\,\vdash_{\class{K},\btau}$ is finitary.
\end{lem}
\noindent
(The finitely algebraizable finitary logics coincide with the original `algebraizable logics'
of Blok and Pigozzi \cite{BP89}.  For cases excluded by their definition, see \cite{Her96,Raf10}.)

Notice that
(\ref{eq:prevariety of logic}) is valid in a class $\class{K}$ of similar algebras iff
it is valid in the prevariety $\mathbb{ISP}(\class{K})$ generated by $\class{K}$.  Our focus on prevarieties is therefore
not restrictive.  In contrast with the case of quasivarieties, it is not provable in the class theory NBG (with choice)
that every prevariety has an  axiomatization involving only a \emph{set} of variables \cite{Adamek1990}.
Papers dealing with the algebraization of logics over proper classes of variables include \cite{BH06,CP99,MRW}.

\begin{remk}\label{not all}
\textup{A
variety $\class{K}$ that
satisfies $f(x,x,\dots,x)\approx x$ for each of its basic operation symbols $f$
is said to be \emph{idempotent}.
In this case, if (\ref{eq:prevariety of logic}) is valid in $\class{K}$, then $\class{K}$ satisfies
${\de_i(x)\approx x\approx \ep_i(x)}$ for all $i\in I$, making the left hand side of (\ref{eq:prevariety of logic})
true on any interpretation of $x,y$ in any member of $\class{K}$.  In view of the right hand side of (\ref{eq:prevariety of logic}),
this forces $\class{K}$ to be trivial.  Thus, no nontrivial idempotent variety is a variety of logic.
In particular, a nontrivial variety of lattices cannot be a variety of logic, as it is idempotent.
More strikingly, although De Morgan lattices are a modest generalization of Boolean algebras, in essentially the same signature,
the variety of De Morgan lattices (which is not idempotent) also fails to be a variety of logic \cite{F97}.
}\qed
\end{remk}




\section{Category Equivalences}

A class $\class{K}$ of similar algebras can be treated as a concrete category, the morphisms being the algebraic homomorphisms
between members of $\class{K}$.  Termwise equivalent classes are then
categorically equivalent,
but not conversely.

When a prevariety $\class{K}$ algebraizes a
logic $\,\vdash$,
we sometimes discover significant features of $\,\vdash$ via `bridge theorems'
of the form
\begin{equation}\label{bridge theorem}
\text{$\vdash$ has metalogical property $P$ iff $\class{K}$ has algebraic property $Q$.}
\end{equation}
Examples include
connections between metalogical interpolation properties
and algebraic amalgamation properties \cite{CP99}, between definability theorems and
the surjectivity of suitable epimorphisms \cite{BH06,MRW}, and between deduction-like theorems and
congruence extensibility properties \cite{BP88,BP,Cze01}.

As it happens, the algebraic properties $Q$ alluded to here are \emph{categorical}, i.e., they persist under
category equivalences between classes $\class{K}$ of the kind to which (\ref{bridge theorem}) applies.
In such cases, if we wish to establish $P$ for $\,\vdash$, we are not forced to prove $Q$ in $\class{K}$ directly;
it suffices to prove $Q$ in an equally suitable class
$\class{M}$ that is \emph{categorically equivalent} to $\class{K}$.\,\footnote{\,The value of this observation
lies not only in the hope that $\class{M}$ can be chosen
simpler or better-understood than $\class{K}$, but also in the possibility that $\class{M}$
algebraizes a logic $\,\vdash'$, different from $\,\vdash$ (perhaps in a different signature).  In that
situation, a category equivalence $F$ between $\class{M}$ and $\class{K}$ carries positive and negative
results from one whole \emph{family} of logics to another.  This is because, in the case of varieties for instance, $F$ induces
an isomorphism between the respective
sub(quasi)variety lattices of $\class{M}$ and $\class{K}$, along which categorical properties can still
be transferred.  And the subquasivarieties of $\class{K}$ [resp.\ $\class{M}$] algebraize the extensions of $\,\vdash$
[resp.\ $\,\vdash'$], with subvarieties corresponding to axiomatic extensions.}

That being so, and in view of Remark~\ref{not all}, it is natural to ask which prevarieties are
categorically equivalent to prevarieties of logic.  We proceed to prove that this is true of \emph{every}
prevariety.
\begin{defn}\label{matrix powers}
\textup{Given an algebra $\sbA$ and $n \in \omega=\{0,1,2,\dots\}$, we denote by $T_{n}(\A)$ the set of all $n$-ary terms in the
signature of $\A$. For $n>0$,
the $n$\textit{-th matrix power} of $\A$ is the algebra
\[
\A^{[n]} \coloneqq \langle A^{n}, \{ m_{t} : t \in T_{kn}(\A)^{n} \text{ for some positive }k \in \omega \},
\]
where for each $t = \langle t_{1}, \dots, t_{n}\rangle \in T_{kn}(\A)^{n}$, we define $m_{t} \colon (A^{n})^{k} \mrig A^{n}$ as follows: if $a_{j} = \langle a_{j1}, \dots, a_{jn}\rangle \in A^{n}$ for $j = 1, \dots, k$, then
\[
m_{t}(a_{1}, \dots, a_{k}) = \langle t^{\A}_{i}(a_{11}, \dots, a_{1n}, \dots, a_{k1}, \dots, a_{kn}) :1 \leq  i \leq n \rangle.
\]
(Roughly speaking, therefore, the basic operations of $\A^{[n]}$ are all conceivable
operations on $n$-tuples that can be defined using the terms of $\sbA$.)}

\textup{For $0<n\in\omega$, the $n$\textit{-th matrix power} of a class $\mathsf{K}$ of similar algebras is the class
$\mathsf{K}^{[n]} \coloneqq \III \{ \A^{[n]} : \A \in \mathsf{K} \}$.}\qed
\end{defn}

Applications of the matrix power construction
in universal algebra range
from the algebraic description of category equivalences and adjunctions \cite{McK96,Mor16b} to the study of clones \cite{Ne70}, Maltsev conditions \cite{Ta75a,GaTa84}, and finite algebras \cite{HM88}.  Matrix powers are also the basis for `twist-product' constructions and product representations; see for instance \cite{CaPri16a}.

\begin{thm}\textup{(cf.\ \cite[Thm.\ 2.3]{McK96})}\label{Thm:equivalence}
Let\/ $\mathsf{K}$ be a class of similar algebras and\/ $n$ a positive integer. Then\/ $\mathsf{K}^{[n]}$ is a class of similar algebras, which is categorically equivalent to\/ $\mathsf{K}$.  Moreover, if\/ $\mathsf{K}$ is a prevariety\/ \textup{[}resp.\ a quasivariety; a variety\/\textup{],} then so is\/ $\mathsf{K}^{[n]}$\textup{.}
\end{thm}

\begin{proof}
It is not difficult to see that the functor $(\cdot)^{[n]} \colon \mathsf{K} \to \mathsf{K}^{[n]}$ sending algebras $\A \in \mathsf{K}$ to $\A^{[n]} \in \mathsf{K}^{[n]}$ and replicating homomorphisms componentwise is a category equivalence.  And for each class operator
$\mathbb{O}$ among $\SSS, \PPP, \PPU, \HHH$, it is easily verified that $\mathsf{K}$ is closed under $\mathbb{O}$ iff the same is true of
$\mathsf{K}^{[n]}$.
\end{proof}

We can now prove the main result of this section.




\begin{thm}\label{Thm:WeakFact}
Let\/ $\class{K}$ be any prevariety.  Then\/ $\class{K}$ is categorically equivalent to a prevariety of logic, i.e., to the
equivalent algebraic semantics\/ $\class{M}$ of some algebraizable logic\/ $\,\vdash$\textup{.}

Moreover, we can choose\/ $\class{M}$
in such a way that the transformers\/ $\btau$ and\/ $\brho$ in\/ $\textup{(\ref{eq:prevariety of logic})}$ are finite, and we
can arrange that\/ $\class{M}$ is a\/ \textup{[}quasi\/\textup{]}variety if\/ $\class{K}$ is.

If\/ $\class{K}$ is a quasivariety,
then\/ $\,\vdash$ can be chosen finitary.
%
\end{thm}
\begin{proof}
Let $\class{M}$ be the matrix power $\class{K}^{[2]}$.  By
Theorem \ref{Thm:equivalence} and Lemma~\ref{Lem:algebraizability}, we need only prove that $\class{M}$
satisfies (\ref{eq:prevariety of logic}) for some finite transformers $\btau,\brho$ (in which case $\,\vdash_{\class{M},\btau}$
can serve as $\,\vdash$).

Now each member of $\class{M}$ has basic binary operations $\to$ and $\leftarrow$, and a basic unary operation $\Box$
such that, for all $\sbA\in\class{K}$ and
$a,b,c,d\in A$,
\begin{align*}
\langle a, b\rangle \to^{\A^{[2]}}\langle c, d \rangle &= \langle a, c \rangle = \langle \pi_{1}(a,b,c,d),
\,\pi_{3}(a,b,c,d)
\rangle;\\
\langle a, b\rangle \leftarrow^{\A^{[2]}}\langle c, d \rangle &= \langle b, d \rangle=\langle \pi_{2}(a,b,c,d),
\,\pi_{4}(a,b,c,d)
\rangle;\\
\Box^{\A^{[2]}}\langle a, b \rangle &= \langle b, a \rangle=\langle \pi_{2}(a,b),\,
\pi_{1}(a,b)
\rangle,
\end{align*}
where $\pi_{k}(z_{1}, \dots, z_{n})\seteq z_k$
whenever $1\leq k\leq n\in\omega$.
These are indeed basic operations for $\class{M}$, because projections are term functions of $\sbA$.

For every $\A \in \class{K}$ and $a, b, c, d \in A$, we have
\begin{eqnarray*}
\langle a, b \rangle = \langle c, d \rangle & \textup{ iff }
& a= c \text{ and }b=d\\
& \textup{ iff }
& \langle a, c\rangle = \langle c, a\rangle \text{ and }\langle b, d\rangle = \langle d, b \rangle\\
& \textup{ iff }
& \big(\langle a, b\rangle \to^{\A^{[2]}} \langle c, d\rangle = \Box^{\A^{[2]}} ( \langle a, b\rangle \to^{\A^{[2]}} \langle c, d\rangle )\text{ and}\\
&  & \,\,\langle a, b\rangle \leftarrow^{\A^{[2]}} \langle c, d\rangle = \Box^{\A^{[2]}} ( \langle a, b\rangle \leftarrow^{\A^{[2]}} \langle c, d\rangle )\big).
\end{eqnarray*}
This implies that the following formula is valid in $\class{M}$:
\begin{equation}
\big(x \to y \approx \Box(x \to y) \;\;\&\;\; x \leftarrow y \approx \Box ( x \leftarrow y)\big) \Longleftrightarrow
x \approx y.
\end{equation}
In other words, (\ref{eq:prevariety of logic}) becomes valid in $\class{M}$ when we set
\[
\btau(x)=\{\langle x,\Box x\rangle\}
\textup { \,and\, }
\brho(x,y)=\{x\to y,\,x\leftarrow y\}.\,\,\footnote{\,Readers who are familiar with
abstract algebraic logic will notice that the reduced matrix models of $\,\vdash_{\class{M},\btau}$ are, up to isomorphism, just all
$\langle\sbA,\{\langle a,a\rangle:a\in A\}\rangle$, $\sbA\in\class{K}$.}
\qedhere
\]
\end{proof}

\begin{cor}
The property of being the equivalent algebraic semantics of an algebraizable logic is not preserved by category equivalences between prevarieties, quasivarieties or varieties.
\end{cor}
\begin{proof}
This follows from Theorem~\ref{Thm:WeakFact} and Remark~\ref{not all}.
\end{proof}



\section{Congruence Equations}

We have noted that the transformer $\brho$ in the definition of a \emph{quasivariety} of logic can be chosen finite.
By a \emph{finitary variety of logic}, we mean a variety of logic for which the transformer $\btau$ can also
be chosen finite (i.e., $\,\vdash_{\class{K},\btau}$ is finitary---see Lemma~\ref{Lem:algebraizability}).
\begin{remk}\label{maltsev class}
\emph{The finitary varieties of logic constitute a Maltsev class in the sense of \cite{Tay73}.  Indeed, suppose
\[
\textup{$\btau=\{\langle \de_i,\ep_i\rangle:i=1,\dots,n\}$ \,and\, $\brho=\{\rho_j:j=1,\dots,m\}$.}
\]
Applying Maltsev's Lemma (cf.\ \cite[Lem.\,V.3.1]{BS81}) to the free $2$-generated algebra in a variety $\class{K}$, we see that
(\ref{eq:prevariety of logic})
is equivalent, over
$\class{K}$, to the conjunction of the identities
$\de_i(\rho_j(x,x))\approx\ep_i(\rho_j(x,x))$ and
a suitable scheme of identities
\begin{eqnarray*}
& x\,\approx\,t_1(x,y,
\ov{\de\rho}
(x,y),\ov{\ep\rho}
(x,y))\\
& t_i(x,y,\ov{\ep\rho}
(x,y),\ov{\de\rho}
(x,y))
\,\approx\,
t_{i+1}(x,y,\ov{\de\rho}
(x,y),\ov{\ep\rho}
(x,y))\quad (1\leq i<k)\\
& t_k(x,y,\ov{\ep\rho}
(x,y),\ov{\de\rho}
(x,y)
)\,\approx\,y
\end{eqnarray*}
involving terms $t_1,\dots,t_k$, where $\ov{\de\rho}
(x,y)$ [resp.\ $\ov{\ep\rho}
(x,y)$] abbreviates
\begin{eqnarray*}
&
\de_1(\rho_1(x,y)),\dots,\de_n(\rho_1(x,y)),\dots,\de_1(\rho_m(x,y)),\dots,\de_n(\rho_m(x,y))\\
&
\textup{[resp.\ $\ep_1(\rho_1(x,y)),\dots,\ep_n(\rho_1(x,y)),\dots,\ep_1(\rho_m(x,y)),\dots,\ep_n(\rho_m(x,y))$].}
\end{eqnarray*}
\indent
If we leave the term symbols $\de_i,\ep_i,\rho_j,t_r$ unspecified, then the finite
conjunction above defines a strong Maltsev class, which need not be idempotent (e.g., the variety
$\class{M}$ in the proof of Theorem~\ref{Thm:WeakFact} does not satisfy $\Box x\approx x$ when $\class{K}$
is a variety).
There are only denumerably many such formal conjunctions,
and any two of them have a common weakening of the same form, got by maximizing,
for each of the letters $\de,\ep,\rho,t$, the number of subscripted occurrences of that letter.
The finitary varieties of logic are therefore directed by the interpretability relation,
whence they form a Maltsev
class.\,\footnote{\,Alternatively, it can be verified that the non-indexed product of two finitary varieties
of logic is a finitary variety of logic; see \cite{Jon80,Neu74,Tay73} for the pertinent characterizations.}}\qed
\end{remk}
Despite this observation, Theorem~\ref{Thm:WeakFact} allows us to show, by an elementary argument, that finitary
varieties of logic are not forced to satisfy any interesting `congruence equation' in the sense of the next
definition.  This contrasts with the fact that the most familiar varieties of logic---the `point-regular'
varieties---are congruence modular and congruence $n$-permutable for a suitable finite $n$ \cite{Hag73}.
(They include the varieties in Example~\ref{boolean etc}.)
\begin{defn}\label{congruence equation}
\textup{A \emph{congruence equation} is a formal equation in the binary symbols $\land$, $\lor$ and $\circ$. It is \emph{satisfied} by an algebra $\A$ if it becomes true whenever we interpret the variables of the equation as congruence relations of $\A$, and for arbitrary binary relations $\alpha$ and $\beta$ on $A$, we interpret $\alpha \land \beta$, $\alpha \lor \beta$ and $\alpha \circ \beta$ as $\alpha \cap \beta$, 
\,$\Theta^{\A}(\alpha \cup \beta)$ and the relational product, respectively.  (Here, $\Theta^\sbA$ stands for congruence generation in $\sbA$.) A congruence equation is \textit{satisfied} by a class of algebras if it is satisfied by every member of the class. It is \textit{nontrivial} if some algebra fails to satisfy it.}\qed
\end{defn}



Because $\circ$ is not generally a binary operation on congruences, we associate with each algebra $\A$ another
algebra
$\Rel(\A) =
\langle \mathit{Rel}(A); \cap,
\lor, \circ \rangle$,
where $\mathit{Rel}(A)$ is the set of all binary relations on $A$, and
\[
\textup{$\alpha \lor \beta \coloneqq
\Theta^\sbA
(\alpha \cup \beta)$ for all $\alpha, \beta \in \mathit{Rel}(A)$.}
\]
The congruence lattice of $\sbA$ is therefore a subalgebra of the $\cap,\lor$ reduct of $\Rel(\A)$.
Given $\alpha,\beta\in\mathit{Rel}(\sbA)$, we also define
\[
\alpha\otimes\beta=\{\langle\langle a,b\rangle,\langle c,d\rangle\rangle:
\langle a,c\rangle\in\alpha \textup{ and }\langle b,d\rangle \in\beta\}\in\mathit{Rel}(A^2).
\]
For \emph{congruences} $\alpha,\beta$ of $\sbA$, it is well known that $\alpha\otimes\beta$ is a congruence
of the algebra $\sbA^2$, but in fact it is also a congruence of $\sbA^{[2]}$.  This follows straightforwardly
from the definitions of $\sbA^{[2]}$ and $\alpha\otimes\beta$.


Recall that a \emph{polynomial}
of an algebra $\langle A;F\rangle$
is a term function of the algebra $\langle A;F\cup F_0\rangle$,
where $F_0$ consists of the elements of $A$, considered as nullary basic operations.
(Of course, we arrange first that $A\cap F=\emptyset$.)

\begin{lem}
Let\/ $\A$ be an algebra.  Then\/ $\lambda \colon \alpha\mapsto\alpha\otimes\alpha$
defines an embedding of\/ $\Rel(\A)$ into\/ $\Rel(\A^{[2]})$\textup{,} which maps congruences to congruences.
%
\end{lem}

\begin{proof}
It is easily verified that, as a function from $\Rel(\A)$ to $\Rel(\A^{[2]})$, $\lambda$ is injective,
$\cap$-preserving and $\circ$-preserving.
Let $\alpha, \beta \in \textup{Rel}(A)$.  We have already mentioned that $\lambda$ preserves congruencehood, from which it follows
that $\lambda(\alpha) \lor \lambda(\beta)\subseteq \lambda(\alpha \lor \beta)$.  It remains to prove the reverse inclusion.


Accordingly, let
$\langle \langle a, b \rangle, \langle c, d \rangle \rangle \in \lambda(\alpha \lor \beta)$, so
$\langle a, c \rangle, \langle b, d \rangle \in \alpha \lor \beta$.
The closure operator $\Theta^\sbA$ (on the power set of $A^2$) is algebraic, so there exist
\begin{equation}\label{e and g}
\langle e_{1}, g_{1}\rangle, \dots, \langle e_{m}, g_{m}\rangle \in \alpha \text{ \,and\, }\langle e_{m+1}, g_{m+1}\rangle, \dots, \langle e_{2m}, g_{2m}\rangle \in \beta
\end{equation}
with
$\langle a, c \rangle, \langle b, d \rangle \in \Theta^\sbA
\{ \langle e_{1}, g_{1}\rangle, \dots, \langle e_{m}, g_{m}\rangle, \langle e_{m+1}, g_{m+1}\rangle, \dots, \langle e_{2m}, g_{2m}\rangle\}$ (where $m$ is finite).
By Maltsev's Lemma, therefore,
there are finitely many $4m$-ary polynomials
$p_{1}, \dots, p_{k}$ and $q_{1}, \dots, q_{k}$ of $\sbA$ such that
\begin{align*}
a &= p_{1}(e_{1}, \dots, e_{2m}, g_{1}, \dots, g_{2m})\\
p_{i}(g_{1}, \dots, g_{2m}, e_{1}, \dots, e_{2m}) &= p_{i+1}(e_{1}, \dots, e_{2m}, g_{1}, \dots, g_{2m})
\\
p_{k}(g_{1}, \dots, g_{2m}, e_{1}, \dots, e_{2m})&= c\,;\\
b &= q_{1}(e_{1}, \dots, e_{2m}, g_{1}, \dots, g_{2m})\\
q_{i}(g_{1}, \dots, g_{2m}, e_{1}, \dots, e_{2m}) &= q_{i+1}(e_{1}, \dots, e_{2m}, g_{1}, \dots, g_{2m})\\
q_{k}(g_{1}, \dots, g_{2m}, e_{1}, \dots, e_{2m})&= d
\end{align*}
for $i = 1, \dots, k-1$.  For each $i\in\{1,\dots,k-1\}$, the rules
%
\begin{align*}
\widehat{p}_{i}(x_{1}, \dots, x_{4m}) \!\coloneqq\, & p_{i}(x_{1}, x_{3}, \dots, x_{2m-1}
, x_{2}, x_{4}, \dots, x_{2m},\\
& x_{2m + 1}, x_{2m + 3}, \dots, x_{4m-1}
, x_{2m+ 2}, x_{2m+4}, \dots, x_{4m});\\[0.25pc]
\widehat{q}_{i}(x_{1}, \dots, x_{4m}) \!\coloneqq\, & q_{i}(x_{1}, x_{3}, \dots, x_{2m-1}
, x_{2}, x_{4}, \dots, x_{2m},\\
& x_{2m + 1}, x_{2m + 3}, \dots, x_{4m-1}
, x_{2m+ 2}, x_{2m+4}, \dots, x_{4m});\\[0.25pc]
t_{i}(z_{1}, \dots, z_{2m}) \coloneqq\! \langle & \widehat{p}_{i}(\pi_1(z_1),\pi_2(z_1),\dots,\pi_1(z_{2m}),\pi_2(z_{2m})),\\
& \widehat{q}_{i}(\pi_1(z_1),\pi_2(z_1),\dots,\pi_1(z_{2m}),\pi_2(z_{2m}))
\rangle
\end{align*}
define two new $4m$-ary polynomials of $\sbA$ and a $2m$-ary polynomial $t_i$ of $\A^{[2]}$,
such that for any
$\langle s_{1}, u_{1}\rangle, \dots, \langle s_{2m}, u_{2m}\rangle \in A^{2}$, the respective first and second co-ordinates
of  $t_{i}
(\langle s_{1}, u_{1}\rangle, \dots, \langle s_{2m}, u_{2m}\rangle)$ are
\begin{align*}
p_{i}
(s_{1}, s_{2}, \dots, s_{m}, u_{1}, \dots, u_{m}, s_{m+1}, \dots, s_{2m}, u_{m+1}, \dots, u_{2m})\,\\
\textup{and\, }q_{i}
(s_{1}, s_{2}, \dots, s_{m}, u_{1}, \dots, u_{m}, s_{m+1}, \dots, s_{2m}, u_{m+1}, \dots, u_{2m}).
\end{align*}
It follows that
\begin{align*}
\langle a, b \rangle =\, & \,t_{1}
(\langle e_{1}, e_{m+1}\rangle, \dots, \langle e_{m}, e_{2m}\rangle, \langle g_{1}, g_{m+1}\rangle, \dots, \langle g_{m}, g_{2m}\rangle);
\end{align*}
\begin{align*}
&
\,t_{i}
(\langle g_{1}, g_{m+1}\rangle, \dots, \langle g_{m}, g_{2m}\rangle, \langle e_{1}, e_{m+1}\rangle, \dots, \langle e_{m}, e_{2m}\rangle)\\
=\, & \,t_{i+1}
(\langle e_{1}, e_{m+1}\rangle, \dots, \langle e_{m}, e_{2m}\rangle, \langle g_{1}, g_{m+1}\rangle, \dots, \langle g_{m}, g_{2m}\rangle);
\end{align*}
\begin{align*}
\,t_{k}
(\langle g_{1}, g_{m+1}\rangle, \dots, \langle g_{m}, g_{2m}\rangle, \langle e_{1}, e_{m+1}\rangle, \dots,
\langle e_{m}, e_{2m}\rangle)
= \langle c, d \rangle
\end{align*}
for $i=1,\dots, k-1$, whence
\begin{align}\label{Eq:cong-gen}
\!\!\!\!\!\!\!\langle \langle a, b \rangle, \langle c, d \rangle \rangle \in
\Theta^{\A^{[2]}}\{ \langle \langle e_{1}, e_{m+1}\rangle, \langle g_{1}, g_{m+1}\rangle \rangle, \dots,  \langle \langle e_{m}, e_{2m}\rangle, \langle g_{m}, g_{2m}\rangle \rangle \}.
\end{align}


Now let $j\in\{1,\dots,m\}$.  By (\ref{e and g}),
\begin{equation}\label{more e and g}
\textup{$\langle \langle e_{j}, e_{j}\rangle, \langle g_{j}, g_{j}\rangle \rangle \in \lambda(\alpha)$ \,and\, $\langle \langle e_{m+j}, e_{m+j}\rangle, \langle g_{m+j}, g_{m+j}\rangle \rangle \in \lambda(\beta)$.}
\end{equation}
For
the basic operation
$\to^{\sbA^{[2]}}$ in the proof of Theorem \ref{Thm:WeakFact}, we have
\begin{align*}
\langle e_{j}, e_{m+j}\rangle &= \langle e_{j}, e_{j}\rangle \to^{\A^{[2]}}
\langle e_{m+j}, e_{m+j}\rangle\\
 \langle g_{j}, g_{m+j}\rangle &= \langle g_{j}, g_{j}\rangle \to^{\A^{[2]}}
 \langle g_{m+j}, g_{m+j}\rangle,
\end{align*}
so
$\langle \langle e_{j}, e_{m+j}\rangle, \langle g_{j}, g_{m+j}\rangle \rangle \in \lambda(\alpha) \lor \lambda(\beta)$,
by (\ref{more e and g}). \vspace{0.7mm}
Then, since $j\in\{
1, \dots, m\}$ was arbitrary,
(\ref{Eq:cong-gen}) yields
$\langle \langle a, b \rangle, \langle c, d \rangle \rangle \in \lambda(\alpha) \lor \lambda(\beta)$, as required.
%
\end{proof}

\begin{cor}\label{Cor:the-trick}
If a congruence equation fails in an algebra $\A$\textup{,} then it fails in $\A^{[2]}$\textup{.}
\end{cor}

This allows us to
prove the main result of this section:

\begin{thm}\label{Thm:congruence-equations}
Every nontrivial congruence equation
fails in some finitary variety of logic, i.e., in a variety that is
the equivalent algebraic semantics of some
finitely algebraizable finitary logic.
\end{thm}
\begin{proof}
Each nontrivial congruence equation
fails in some variety $\mathsf{K}$, hence also
in $\mathsf{K}^{[2]}$ (by Corollary \ref{Cor:the-trick}), which is itself a
variety (by Theorem~\ref{Thm:equivalence}).
And $\mathsf{K}^{[2]}$ is a finitary variety of logic, by
the proof
of Theorem \ref{Thm:WeakFact}.
\end{proof}

A finitary variety of logic satisfying no nontrivial congruence equation in the signature $\land,\lor$ (excluding $\circ$)
was exhibited in \cite{BR01}.  The stronger fact that this variety satisfies no nontrivial idempotent Maltsev condition
was pointed out in \cite[Sec.~10.1]{BR08}, using \cite[Thm.~4.23]{KK}.  Theorem~\ref{Thm:congruence-equations} does not follow
from these observations and general results of universal algebra, however, because it is not evident that every nontrivial congruence
equation (in the full signature $\land,\lor,\circ$) entails a nontrivial idempotent Maltsev condition, as opposed to a weak
Maltsev condition.  (This corrects an impression left in the last lines of \cite[p.~647]{BR08},
and in \cite{Raf07}.)  For more on the general connections between these notions, see \cite{KK}.

\bibliographystyle{plain}

\end{document}